\documentclass{amsart}
\newtheorem{thm}{Theorem}[section]
\newtheorem{lem}{Lemma}[section]
\newtheorem{cor}{Corollary}[section]

\theoremstyle{remark}
\newtheorem{rem}{Remark}[section]
\theoremstyle{definition}

\def\NN{{\Bbb N}}

\def\XX{{\Bbb X}}
\def\CC{{\Bbb C}}

\def\Bf{\mathcal{B}}
\def\Df{\mathcal{D}}
\def\Pf{\mathcal{P}}

\def\Ef{\mathcal{E}}
\def\Ff{\mathcal{F}}

\renewcommand{\Re}{{\Bbb R}}

\def\1{1\!\!\hbox{{\rm I}}}

\def\eqdef{\mathop{=}\limits^{df}}
\def\ax{\Re^+}
\newcommand{\demo}{\emph{Proof.} }
\newcommand{\prt}{\partial}
\newcommand{\be}{\begin{equation}}
\newcommand{\ee}{\end{equation}}
\newcommand{\ba}{\begin{aligned}}
\newcommand{\ea}{\end{aligned}}

\begin{document}

\title[Poincar\'e inequality and exponential integrability]
    {Poincar\'e inequality and exponential integrability of the  hitting times of a Markov process}

\author{Alexey M. Kulik}
\address{Kiev 01601 Tereshchenkivska str. 3, Institute of Mathematics,
Ukrai\-ni\-an National Academy of Sciences}
\email{kulik @imath.kiev.ua}


\begin{abstract}
Extending the  approach of the paper [Mathieu, P. (1997) Hitting times and
spectral gap inequalities, Ann. Inst. Henri Poincar\'e 33, 4, 437
-- 465], we prove that the Poincar\'e inequality for a (possibly non-symmetric) Markov process yields the exponential integrability of the  hitting times of this process. For  symmetric elliptic diffusions, this provides a criterion for the Poincar\'e inequality in the terms of  hitting times.
\end{abstract}

\maketitle

\section{Introduction}
In this paper, we investigate the relations between the following two topics:

\begin{itemize}
\item  the \emph{Poincar\'e inequality} for the
 Dirichlet form associated to a Markov process;

\item the exponential integrability for \emph{hitting times} of the process.
\end{itemize}

It is well known  that when the state space of a (symmetric) Markov process is finite, the topics listed above are, in fact, equivalent; see the
detailed exposition in \cite{AF}, Chapters 2 -- 4. When the state space is infinite, the relations between these topics are more complicated, see the detailed discussion in \cite{Kul11}. The purpose of this paper is two-fold. First, we extend the approach of the paper
\cite{Mat97}, where it was proved that hitting times of a Markov process are integrable assuming some weak version of the Poincar\'e inequality holds true. In this paper, we prove that the Poincar\'e inequality itself provides the \emph{exponential integrability} for hitting times. Next, we show that, on the contrary, under some additional assumptions the  {exponential integrability} for hitting times yields the Poincar\'e inequality. According to the recent paper \cite{Kul11}, for a given Markov process $X$ the \emph{spectral gap}
property  can be verified in the following way. First, one proves some of the \emph{local non-degeneracy} conditions on the transition  probabilities of the Markov process $X$ (minorization condition, Doeblin condition, Dobrushin condition).
Second, one finds {\it
some} Lyapunov-type function $\phi$ such that the recurrence conditions 1) -- 3) of
Theorem 2.2 in \cite{Kul11} hold true. Then $X$ admits an \emph{exponential $\phi$-coupling} (see Definition 2.2 in \cite{Kul11}) according to Theorem 2.2 in \cite{Kul11}. When $X$ is time-reversible (i.e. symmetric), this provides the spectral gap property, which in this case is equivalent to the the Poincar\'e inequality.  The case where $X$ is time-irreversible is more intrinsic; because we do not address this case in the current paper we skip the discussion of this case and refer the interested reader to the paper \cite{Kul11}. Note that in the
strategy  outlined above there is a lot of freedom in the choice of the Lyapunov-type function
$\phi$.  Proposition 2.4 in \cite{Kul11} reduces such a choice  to the
class of the functions of the form
\be\label{01}\phi(x)=E_xe^{\alpha\tau_K},\quad  \hbox{where} \quad
\tau_K=\inf\{t: X_t\in K\}
\ee
is the hitting time of a compact set $K$. In Section 3 below we demonstrate that, for particularly important class of symmetric diffusions,
this  reduction
corresponds to the matter of the problem precisely, and deduce the Poincar\'e inequality under the  exponential integrability assumption for the hitting times. Together with the result of Section 2 this provides an equivalence for symmetric diffusions of the  Poincar\'e inequality on one hand, and the the  exponential integrability of the hitting times on the other hand. Such an equivalence for linear diffusions was established in the recent preprint \cite{Louk09}; note that the case of the symmetric elliptic diffusion on a multidimensional manifold, considered in the current paper, is more complicated because one can not apply here such particularly useful characteristics  of the linear diffusion as the scale function and the speed measure. It should be mentioned that in the recent years the relations between the functional inequalities related to a Markov process (like the Poincar\'e inequality) and the ergodic properties of this process (like the Lyapunov-type condition)  have been studied extensively, e.g. \cite{BCG08}, \cite{CGWW09}, \cite{CGA10}. In this paper, we show, in particular, that for a symmetric elliptic diffusion this relation is one-to-one, and the Poincar\'e inequality is \emph{equivalent} to the Lyapunov-type condition with the Lyapunov-type function of the form (\ref{01}) (this equivalence was investigated, as well, in the recent preprint \cite{CGA10}).

\section{Exponential moments for hitting times under Poincar\'e inequality}

We consider a time homogeneous Markov process $X=\{X_t,
t\in\ax\}$ with a locally compact metric space $(\XX,\rho)$ as
the state space. The process $X$ is supposed to be strong Markov
and to have c\'adl\'ag trajectories. The transition function for
the process $X$ is denoted by $P_t(x,dy), t\in\ax, x\in\XX$. We
use standard notation $P_x$ for the distribution of the process
$X$ conditioned that $X_0=x, x\in\XX$, and $E_x$ for the expectation w.r.t.
$P_x$. All the functions on $\XX$ considered in the paper are assumed to
be measurable w.r.t.  {the} Borel $\sigma$-algebra $\Bf(\XX)$. The set
of probability measures on $(\XX,\Bf(\XX))$ is denoted by
$\Pf(\XX)$. For a given $\mu\in\Pf(\XX)$ and $t\in\ax$, we denote
$ \mu_t(dy)\eqdef \int_\XX P_t(x,dy)\, \mu(dx). $
 {The} probability  measure $\mu$ is called an invariant measure for $X$
if  $\mu_t=\mu, t\in\ax$.

In what follows, we suppose an invariant measure $\pi$ for the
process $X$ to be fixed. The process $X$ generates the semigroup $\{T_t\}$ in $L_2=L_2(\XX, \pi)$:
$$
T_tf(x)=\int_\XX f(y)P_t(x,dy), \quad f\in L_2, \quad t\in \ax.
$$
Let $A$ be the generator of this semigroup, and $\Ef$ be the associated
Dirichlet form; which is defined  as  the completion of the  bilinear form
$$
Dom(A)\times Dom(A)\ni (f,g)\mapsto -(Af,g)_{L_2}
$$
with respect to the norm $\|\cdot\|_{\Ef,1}\eqdef
\Big[\|\cdot\|_{2}^2-(A\cdot,\cdot)_{L_2}\Big]^{1/2}$ (e.g.
\cite{MR92}, Chapter 2). Within this paper, we are mainly interested in the following {\it Poincar\'e inequality}:
\be\label{34} \mathrm{Var}_\pi(f)\eqdef\int_\XX f^2d\pi-\left(\int_\XX f
d\pi\right)^2\leq c \,\Ef(f,f),\quad f\in Dom(\Ef).\ee

In what follows, the form $\Ef$ is supposed to be regular; that is,
the set $Dom(\Ef)\cap C_0(\XX)$ is claimed to be dense both in
 $Dom(\Ef)$ w.r.t. the norm $\|\cdot\|_{\Ef,1}$ and in $C_0(\XX)$
 w.r.t. uniform convergence on a compacts ($C_0(\XX)$ is the set of continuous functions with compact
 supports).  We also assume that the following {\it sector condition} holds true:
 $$
\exists D\in\ax:\quad  |\Ef( f,g)|\leq D
\|f\|_{\Ef,1}\|g\|_{\Ef,1},\quad f,g\in Dom(\Ef).
 $$

 It is well-known (see the discussion in
Introduction to \cite{Mat97} and references therein) that the
hitting times  $\tau_K$ have natural application in the
probabilistic representation for the family of {\it
$\alpha$-potentials} for the Dirichlet form $\Ef$. The
$\alpha$-potential, for given $\alpha>0$ and closed $K\subset
\XX$, is defined as the function $h_\alpha^K\in Dom(\Ef)$ such
that $h_\alpha^K=1$ quasi-everywhere on $K$, and
$\Ef(h_\alpha^K,u)=-\alpha(h_\alpha^K,u)$ for every
quasi-continuous function $u\in Dom(\Ef)$ such that $u=0$
quasi-everywhere on $K$. On the other hand,
$$
h_\alpha^K(x)=E_xe^{-\alpha \tau_K}, \quad x\in\XX.
$$

It is a straightforward corollary of the part (i) of the main
theorem from \cite{Mat97} that, if $X$ possesses (\ref{34}) with
some $c>0$, then $E_\pi \tau_K<+\infty$  for every $K$ with
$\pi(K)>0$  (here and below, $E_\pi\eqdef \int_\XX
E_x\,\pi(dx)$). We will prove the following  stronger version of
this statement.

\begin{thm}\label{t51} Assume $X$ possess (\ref{34}) with some
$c>0$. Then for every closed set $K\subset \XX$ with
$\pi(K)>0$
$$
E_\pi e^{\alpha\tau_K}<+\infty,\quad  \alpha<{\pi(K)\over
c}.
$$
Moreover, the function $h_{-\alpha}^K(x)\eqdef E_x
e^{\alpha\tau_K}, x\in\XX$ possesses the following properties:

a) $h_{-\alpha}^K\in Dom(\Ef)$ and $h_{-\alpha}^K=1$ on $K$;

b) $\Ef(h_{-\alpha}^K,u)=\alpha(h_{-\alpha}^K,u)$ for every
quasi-continuous function $u\in Dom(\Ef)$ such that $u=0$
quasi-everywhere on $K$.
\end{thm}

 \demo We assume $K$ to be fixed and omit the
 respective index in the notation, e.g. write $\tau$ for $\tau^K$ and $h_\alpha$ for $h_\alpha^K$.  For $z\in \CC$ with
 $\mathrm{Re}\,z>0$, define respective $z$-potential:
 $$
 h_z(x)=E_xe^{-z \tau}, \quad x\in\XX.
 $$

Denote by $H_\Ef$ the  $Dom(\Ef)$ considered as a Hilbert space
with the scalar product $$(f,g)_{\Ef,1}\eqdef
(f,g)_{L_2}+\Ef(f,g).$$ The following lemma shows  that
$\{h_z,\mathrm{Re}\, z>0\}$  can be considered as an analytical
extension of the family of $\alpha$-potentials
$\{h_\alpha,\alpha>0\}\subset H_\Ef$ that, in addition,  keeps
the properties of this family.

\begin{lem}\label{l51} 1) The function $z\mapsto h_z$ is analytic as a function
taking values in the Hilbert space $H_\Ef$.

2) For every $z$ with $\mathrm{Re}\, z>0$, the following
properties hold:

(i) $h_z=1$ quasi-everywhere on $K$;

(ii) $\Ef(h_z,u)=-z(h_z,u)$ for every quasi-continuous function
$u\in Dom(\Ef)$ such that $u=0$ quasi-everywhere on $K$.
\end{lem}

\demo
 Denote $h_z^m(x)=(-1)^m  E_x\tau^m e^{-z \tau}, x\in\XX, m\geq 1$. One can verify easily that, for every $m\in\NN$,
 \be\label{hm}
 {d^m\over dz^m} h_z=h_z^m
\ee
on the set $\{z: \mathrm{Re}\,z>0\}$, with the function $z\mapsto
h_z$ is considered as a function taking values in $L_2$. In
addition,
$$
\|h_z^m\|^2_2\leq E_\pi \left|\tau^m e^{-z \tau}\right|^2=E_\pi
\tau^{2m}e^{-2\tau\mathrm{Re}\, z}\leq {(2m)!\over
(2\mathrm{Re}\,z)^{2m}},
$$
since ${(2\tau\mathrm{Re}\, z)^{2m}\over (2m)!}\leq
e^{2\tau\mathrm{Re}\, z}$. Therefore, \be\label{52}
{\|h^m_z\|_2\over m!}\leq \sqrt{C_{2m}^m\over
2^{2m}}(\mathrm{Re}\, z)^{-m}<(\mathrm{Re}\, z)^{-m},\quad m\in
\NN, \ee and hence the function
$$
\{z: \mathrm{Re}\,z>0\}\ni z\mapsto h_z\in L_2
$$
is analytic.

 For every
$\alpha,\alpha'>0$ we have $h_\alpha-h_{\alpha'}=0$
quasi-everywhere on $K$. Hence \be\label{51}\ba
\Ef(h_\alpha-h_{\alpha'},h_\alpha-h_{\alpha'})&=\Ef(h_\alpha,h_\alpha-h_{\alpha'})-\Ef(h_{\alpha'},h_\alpha-h_{\alpha'})\\
&=-\alpha(h_\alpha,h_\alpha-h_{\alpha'})+\alpha'(h_{\alpha'},h_\alpha-h_{\alpha'})\\
&=(\alpha'-\alpha)(h_{\alpha'},h_\alpha-h_{\alpha'})+ \alpha
(h_{\alpha'}-h_{\alpha},h_\alpha-h_{\alpha'}). \ea \ee

For a given $\alpha>0$ and $\alpha'\to \alpha$, the family
$\{{h_{\alpha'}-h_{\alpha}\over \alpha'-\alpha}\}$ converges to $h_\alpha^1$ in $L_2$, see (\ref{hm}). Then (\ref{51}) yields that this family
is bounded in
$H_\Ef$, and thus is weakly compact in $H_\Ef$. Combined with the fact that this family is converges in $L_2$, this yields
that the function  $(0,+\infty)\in\alpha\mapsto h_\alpha\in
H_\Ef$ is differentiable in a weak sense, and  $h^1_\alpha$
equals its (weak) derivative at the point $\alpha$.

We have $h_\alpha^1=0$ quasi-everywhere on $K$, since
$$
h_\alpha(x)=1\Leftrightarrow e^{-\alpha \tau}=1 \
P_x-\hbox{a.s.}\Leftrightarrow  \tau=0 \
P_x-\hbox{a.s.}\Leftrightarrow h_\alpha^1(x)=0.
$$
In addition, since $h_\alpha^1$ is a weak derivative of
$h_\alpha$, we have
$\Ef(h_\alpha^1,u)=-(h_\alpha,u)-\alpha(h_\alpha^1,u)$ for every
quasi-continuous function $u\in Dom(\Ef)$ such that $u=0$
quasi-everywhere on $K$.  Now, repeating the same arguments, we
get by induction that, for every $m\geq 1$, the function
$(0,+\infty)\in\alpha\mapsto h_\alpha\in H_\Ef$ is $m$ times
weakly differentiable, $h_\alpha^m$ is the corresponding weak
derivative of the $m$-th order, and the following properties hold:

(i$^m$) $h_\alpha^m=0$ quasi-everywhere on $K$;

(ii$^m$)
$\Ef(h_\alpha^m,u)=-(h_\alpha^{m-1},u)-\alpha(h_\alpha^m,u)$ for
every quasi-continuous function $u\in Dom(\Ef)$ such that $u=0$
quasi-everywhere on $K$.

Property (ii) with $u=h_\alpha^m$ and estimate (\ref{52}) yield
that, for a given $\alpha$, series
$$
H_z\eqdef h_\alpha+\sum_{m=1}^\infty{z^m\over m!} h_\alpha^m\in
H_\Ef
$$
converge in the circle $\{|z-\alpha|<\alpha\}$. The sum is a
weakly analytic $H_\Ef$-valued function, and hence is  analytic
(\cite{Rud73}, Theorem 3.31). On the other hand, the same series
converge in $L_2$ to $h_z$. This yields that $h_z=H_z$ in the
circle $\{|z-\alpha|<\alpha\}$. By taking various
$\alpha\in(0,+\infty)$, we get that the function $z\mapsto h_z$
is  an $H_\Ef$-valued analytic function inside the angle
$\Df_1\eqdef\{z:\mathrm{Re}\, z> |\mathrm{Im}\,z|\}$. In
addition, properties (i$^m$), (ii$^m$) of the $m$-th coefficients
of the series ($m\geq 1$) provide that $h_z$ satisfy (i),(ii)
inside the angle.

Now, we complete the proof using the following iterative
procedure. Assume that the function $z\mapsto h_z\in H_\Ef$ is
analytic in some domain $\Df\subset \{z:\mathrm{Re}\,z >0\}$ and
satisfy (i),(ii) in this domain. Then the same arguments with
those used above show that, for every $z_0\in \Df$, the domain
$\Df$ can be extended to
$\Df'\eqdef\Df\cup\{z:|z-z_0|<\mathrm{Re}\,z_0\}$ with the
function $z\mapsto h_z$ still being analytic in $\Df'$ and
satisfying (i),(ii) in the extended domain. Therefore, we prove
iteratively that the required statement holds true in every angle
$\Df_k\eqdef\{z:\mathrm{Re}\, z> {1\over k} |\mathrm{Im}\,z|\}$.
Since $\cup_k\Df_k=\{z:\mathrm{Re}\,z>0\}$, this completes the
proof.

 Next, we consider ``$\psi$-potentials'' that correspond to   functions $\psi:\ax\to \Re$.
Denote
$$h_\psi(x)=E_x\psi(\tau),\quad x\in\XX.
$$
The following statement is an appropriate modification of the
inversion formula for the Laplace transform.

\begin{lem}\label{l52} Let $\psi\in C^2(\Re)$ have a compact support and $\mathrm{supp}\,\psi\subset [0,+\infty)$.
Denote  $\Psi(z)=\int_\Re e^{zt}\psi(t)\, dt,$ $z\in\CC$.

 The function
$h_\psi$ belongs to $H_\Ef$ and admits integral representation
\be\label{53} h_\psi={1\over 2\pi i
}\int_{\sigma-i\infty}^{\sigma+i\infty}\Psi(z)h_z\,dz, \ee where
$\sigma>0$ is arbitrary, and the integral is well defined as an
improper Bochner integral of an $H_\Ef$-valued function.
\end{lem}

\demo First, let us show that the integral in the right hand side
of (\ref{53}) is well defined. We have by condition (ii) of Lemma
\ref{l51} that
$$
\Ef(h_z,h_z)=\Ef(h_z,h_z-1)=-z(h_z,h_z-1).
$$
For any $z$ with $\mathrm{Re}\,z>0$, we have $|h_z(x)|\leq E_x
e^{-\tau\mathrm{Re}\,z}\leq 1$, and thus $|h_z(x)-1|\leq 2$.
Hence,
$$
\|h_z\|_{H_\Ef}= \sqrt{\|h_z\|_2^2+\Ef(h_z,h_z)}\leq
\sqrt{1+2|z|}.
$$
On the other hand, for $\psi$ satisfying conditions of the lemma,
$$
z^2\Psi(z)=\int_0^\infty e^{zt}\psi''(t)\, dt,\quad
|z^2\Psi(z)|\leq \int_0^\infty e^{t\mathrm{Re}\,z}|\psi''(t)|\,
dt.
$$
Thus, on the line $\sigma+i\Re\eqdef\{z:\mathrm{Re}\, z=\sigma\}$,
the function $z\mapsto \Psi(z)h_z\in H_\Ef$ admits the following
estimate:
$$
\|\Psi(z)h_z\|_{H_\Ef}\leq C|z|^{-{3\over 2}},
$$
and therefore it is integrable on $\sigma+i\Re$. Denote by
$g_\psi\in H_\Ef$ corresponding integral. In order to prove that
$h_\psi=g_\psi$, it is sufficient to prove that $h_\psi$ and
$g_\psi$ coincide as elements of $L_2$. Hence, we have reduced the
proof of the lemma to verification of the following ``weak
$L_2$-version'' of (\ref{53}): \be\label{54} \int_\XX h_\psi v\,
d\pi={1\over 2\pi i}\int_\XX
\int_{\sigma-i\infty}^{\sigma+i\infty}\Psi(z)h_z(x)v(x)\,dz\pi(dx),\quad
v\in L_2. \ee Recall that $h_z(x)=E_xe^{-z\tau}$, and hence the
right hand side of (\ref{54}) can be rewritten to the form
$$
{1\over 2\pi i}\int_\XX
\int_{\sigma-i\infty}^{\sigma+i\infty}E_x\Psi(z)
e^{-z\tau}v(x)\,dz\pi(dx)={1\over 2\pi i}\int_\XX E_x
\int_{\sigma-i\infty}^{\sigma+i\infty}\Psi(z)
e^{-z\tau}v(x)\,dz\pi(dx).
$$
Here, we have changed the order of integration using Fubini's
theorem. This can be done, because $|\Psi(z)|\leq C|z|^{-2}$, and
therefore
$$
E_x \int_{\sigma-i\infty}^{\sigma+i\infty}|\Psi(z) e^{-z\tau}|\,
dz =h_\sigma(x) \int_{\sigma-i\infty}^{\sigma+i\infty}|\Psi(z)|\,
dz\leq C h_\sigma(x).
$$
The function $\Psi$ is the (two-sided) Laplace transform for
$\psi$, up to the change of variables $p\mapsto -z$. We write the
inversion formula for the Laplace transform in the terms of
$\Psi$ and, after the change of variables, get
$$
\psi(t)={1\over 2\pi
i}\int_{-\sigma-i\infty}^{-\sigma+i\infty}e^{pt}\Psi(-p)\,
dp={1\over 2\pi
i}\int_{\sigma-i\infty}^{\sigma+i\infty}e^{-zt}\Psi(z)\, dz,
\quad t\in \ax.
$$
Hence, the right hand side of (\ref{54}) is equal
$$
\int_\XX E_x\psi(\tau)v(x)\pi(dx)=\int_\XX h_\psi v\, d\pi,
$$
that proves (\ref{54}).

\begin{cor}\label{c51}  Let $\psi\in C^3(\Re)$ and $\mathrm{supp}\,\psi'\subset [0,+\infty)$.
Then $h_\psi\in Dom(\Ef)$ and   \be\label{55}
\Ef(h_\psi,u)=(h_{\psi'},u) \ee for every $u\in Dom(\Ef)$ such
that $u=0$ quasi-everywhere on $K$.
\end{cor}
\demo Assume first that $\int_{\ax} \psi'(x)\, dx=0$.  Then both
$\psi$ and $\psi'$ satisy conditions of Lemma \ref{l52}. We have
$\tilde \Psi(z)\eqdef\int_\Re e^{zt}\psi'(t)\, dt=-z\Psi(z)$.
Hence, from the representation (\ref{53}) for $h_\psi$ and
$h_{\psi'}$ and relation $\Ef(h_z,u)=-z(h_z,u),\mathrm{Re}\,
z>0$, we get
$$
\Ef(h_\psi,u)={1\over 2\pi i
}\int_{\sigma-i\infty}^{\sigma+i\infty}\Psi(z)\Ef(h_z,u)\,dz={1\over
2\pi i }\int_{\sigma-i\infty}^{\sigma+i\infty}\tilde
\Psi(z)(h_z,u)\,dz=(h_{\psi'},u).
$$

The general case can be reduced to the one considered above by
the following limit procedure. Since $\mathrm{supp}\,\psi'\subset
[0,+\infty)$, there exist $C\in \Re$ and  $x_*\in \ax$ such that
$\psi(x)=C, x\geq x_*$. Take a function $\vartheta\in C^3(\Re)$ such that $\vartheta(x)=0, x\leq 0, \vartheta(x)=C, x\geq 1,$
 and put
$$
\psi_t(x)=\psi(x)-\vartheta(x-t),\quad x\in\Re, t>x_*.
$$
Then every $\psi_t$ satisfies the additional assumption
$\int_{\ax} [\psi_t]'(x)\, dx=0$, and thus $h_{\psi_t}$ belongs
to $Dom(\Ef)$ and satisfies (\ref{55}). It can be verified easily
that $h_{\psi_t}\to h_\psi, t\to \infty$ in $L_2$ sense. In
addition,
$$
\Ef(h_{\psi_t},h_{\psi_t})=(h_{[\psi_t]'},h_{\psi_t})\to
(h_{[\psi]'},h_{\psi})<+\infty,\quad t\to+\infty
$$
(here, we have used (\ref{55}) with $u=h_{\psi_t}$). This means
that the family $\{h_{\psi_t}\}$ is bounded in $H_\Ef$, and hence
is weakly compact in $H_{\Ef}$. Therefore, $h_{\psi_t}\to h_\psi,
t\to \infty$ weakly in $H_{\Ef}$. Since $h_{[\psi_t]'}\to
h_{\psi'}, t\to \infty$ in $L_2$ sense, (\ref{55}) for $\psi$
follows from (\ref{55}) for $\psi_t$.

 Now, we are ready to complete the proof of the theorem. Let us  fix
 $\alpha<{\pi(K)\over c}$, and construct the family of the functions $\varrho_t, t\geq 1$ that approximate the function
 $\varrho: x\mapsto e^{\alpha x}-1$ appropriately. First, we take
 function
 $\chi\in C^3(\Re)$ such that $\chi\geq 0, \chi'\leq 0, \chi(x)=1,
 x\leq0$, and $\chi(x)=0, x\geq 1$. We put
 $$
 \rho_t(x)=\int_0^x\alpha e^{\alpha y}\chi(y-t)\, dy, \quad
 x\geq 0, t\geq 1.
 $$
 By the construction, the derivatives of the functions $\rho_t, t\geq
 1$ have the following properties:

 a) $[\rho_t]'\geq 0$ and $[\rho_t]'(x)=0, x\geq t+1$;

 b) $[\rho_s]'\leq[\rho_t]', s\leq t$.

 Since $\rho_t(0)=0, t\geq 1$, the latter property yields that $\rho_s\leq \rho_t, s\leq
 t$. In addition,
 $$
 [\rho_t]''(x)= \alpha e^{\alpha x} \chi'(x)+\alpha^2 e^{\alpha x} \chi(x)\leq
 \alpha^2 e^{\alpha x} \chi(x)=\alpha [\rho_t]'(x),
 $$
 since $\chi'\leq 0$. This and relation $[\rho_t]'(0)=\alpha (\rho_t(0)+1)$
 provide
 \be\label{56} [\rho_t]'\leq \alpha(\rho_t+1).
 \ee

At last,  we take  function $\theta\in C^3(\Re)$ such that
$\theta'\geq 0, \theta(x)=0,
 x\leq 0$, and $\theta(x)=1, x\geq 1$. We put
 $$
 \varrho_t(x)=\begin{cases} \theta\left({xt}\right)\rho_t(x),&x\geq
 0\\
  0,& x<0
 \end{cases},\quad t\geq 1.
 $$
We have $\varrho_t\uparrow \varrho, t\uparrow \infty$. In
addition, by (\ref{56}), \be\label{57}
[\varrho_t]'(x)=t\theta'(tx)\rho_t(x)+\theta(tx)[\rho_t]'(x)\leq
t\sup_{y}\theta'(y)\rho_t(t^{-1})+\alpha(\rho_t(x)+1)\leq
\alpha\varrho_t(x)+C \ee with an appropriate constant $C$ (recall
that $t\rho_t(t^{-1})=t\alpha(e^{\alpha t^{-1}}-1)\to \alpha^2,
t\to \infty$).

Every $\varrho_t$ satisfies conditions of Corollary \ref{c51},
and hence $$
\int_{\XX}h_{\varrho_t}^2\,d\pi-\left(\int_{\XX}h_{\varrho_t}\,d\pi\right)^2\leq
c\Ef(h_{\varrho_t}, h_{\varrho_t})=c (h_{[\varrho_t]'}, h_{\varrho_t})\leq \alpha c(h_{\varrho_t}, h_{\varrho_t})+C\int_\XX h_{\varrho_t}\,
d\pi.$$  Here, we have used subsequently property (\ref{34}),
equality  (\ref{56}) with $u=h_{\varrho_t}$, and (\ref{57}).

We have $h_{\varrho_t}=0$ on $K$ because $\varrho_t(0)=0$. Then,
by the Cauchy inequality,
$$\ba
\int_{\XX}h_{\varrho_t}^2\,d\pi-&\left(\int_{\XX}h_{\varrho_t}\,d\pi\right)^2=
\int_{\XX}h_{\varrho_t}^2\,d\pi-\left(\int_{\XX\setminus K
}h_{\varrho_t}\,d\pi\right)^2\\
&\geq (1-\pi(\XX\setminus K))
\int_{\XX}h_{\varrho_t}^2\,d\pi=\pi(K)(h_{\varrho_t},
h_{\varrho_t}).\ea
$$
Therefore,  $$(h_{\varrho_t}, h_{\varrho_t})\leq {\alpha c \over
\pi(K)}(h_{\varrho_t}, h_{\varrho_t})+C\int_\XX
h_{\varrho_t}\, d\pi,$$ which implies that \be\label{58}
(h_{\varrho_t}, h_{\varrho_t})\leq {C\pi(K)\over \pi(K)-\alpha c}\int_\XX h_{\varrho_t}\, d\pi
\ee (recall that $\alpha<{\pi(K)\over c}$). One can verify
easily that (\ref{58}) yields that the $L_2$-norms of the
functions $h_{\varrho_t}$ are uniformly bounded. Since
$\varrho_t\uparrow \varrho$, this implies that the function
$$
h_\varrho (x)\eqdef E_xe^{\alpha \tau}-1,\quad x\in\XX
$$
belongs to $L_2$, and $h_{\varrho_t}\to h_\varrho,
 t\to \infty$ in $L_2$. Similarly to the proof of Corollary
\ref{c51}, one can verify that  $\{h_{\varrho_t}\}$ is a
bounded subset in $H_\Ef$, and hence $h_{\varrho_t}\to h_\varrho,
 t\to \infty$ weakly in $H_\Ef$. This proves statement a) of the
 theorem.  In order to prove statement b), we apply (\ref{55}) to
$\psi=\varrho_t$, and pass to the limit as $t\to +\infty$. The
 theorem is proved.

\section{Poincar\'e inequality for symmetric diffusions: criterion in the terms of  hitting times}

Let $\XX$ be a connected locally compact  Riemannian manifold of
dimension $d$, and $X$ be a diffusion process on $\XX$. Let $\pi\in \Pf(\XX)$ be an invariant measure for the process $X$
(we assume invariant measure to exist). We assume that $X$ is symmetric w.r.t. $\pi$; that is, $T_t=T_t^*, t\in \ax$.

On a given
local chart of the manifold $\XX$, the generator of the process
$X$ has the form
$$
A=\sum_{j=1}^d a_j\prt_j+{1\over 2}\sum_{j,k=1}^d
b_{jk}\prt^2_{jk},
$$
where $a=\{a_j\}_{j=1}^d$ and $b=\{b_{jk}\}_{j,k=1}^d$ are the
drift and diffusion coefficients of the process $X$ on this chart,
respectively. We assume the coefficients $a,b$ to be H\"older
continuous on every local chart, and the drift $b$ coefficient to
satisfy ellipticity condition
$$
\sum_{j,k=1}^d b_{jk} v_jv_k\geq \beta \sum_{j=1}^dv_j^2
$$
uniformly on every compact. Under these conditions, the
transition function of the process $X$ has a positive density
w.r.t. Riemannian volume, and this density  is a continuous
function on $(0,+\infty)\times \XX\times \XX$. One can easily
deduce this from the same statement for diffusions in $\Re^d$
(e.g. \cite{IKO62}) and strong Markov property of $X$. This
implies that $X$ satisfies  the \emph{extended Doeblin condition} (see Section 2.1 in \cite{Kul11}) on
every compact subset of $\XX$.

\begin{thm}\label{t61} The following statements are equivalent:

1) the Poincar\'e inequality (\ref{34}) holds true with some  constant $c$;

2) the process $X$ admits an exponential $\phi$-coupling for some
function $\phi$, see Definition 2.2 in \cite{Kul11};

3) for every closed subset $K\subset \XX$ with $\pi(K)>0,$ there
exists $\alpha>0$ such that
$$E_\pi e^{\alpha\tau_K}<+\infty.
$$

\noindent In addition, 1) -- 3) hold true assuming that

3$\,'$) there exists a compact subset $K\subset \XX$ and
$\alpha>0$ such that
$$E_x e^{\alpha\tau_K}<+\infty\quad \hbox{for $\pi$-almost all}\quad x\in X.
$$

\end{thm}

\begin{rem} Note that the property 2) both provides  the uniqueness of the invariant measure and makes it possible to give  explicit bounds for the   convergence rate of the transition probabilities of the process $X$ to the invariant distribution, see \cite{Kul11}. Hence for symmetric diffusions the above theorem, together with the criterion for the Poincar\'e inequality in the terms of  hitting times, gives a sufficient condition for an (exponential) ergodicity.
\end{rem}
\demo Implication 2) $\Rightarrow$ 1) follows immediately from Theorem 3.4 in \cite{Kul11}. Implication 1) $\Rightarrow$ 3)
is provided by  Theorem \ref{t51}. Implication 3) $\Rightarrow$ 3$\,'$) is trivial. To prove implication 3$\,'$) $\Rightarrow$ 2),
we will use  Proposition 2.4 in \cite{Kul11}. Recall that  we have
already seen that $X$ satisfies the extended Doeblin condition on
$K$. Hence, in order to apply Proposition 2.4 in \cite{Kul11}, it is sufficient  to verify for a given $\tilde \alpha\in (0,\alpha)$ the following conditions:
\begin{itemize}\item[(a)]  $E_xe^{\tilde \alpha\tau_K}<+\infty,x\in\XX;$
\item[(b)] there exists $S>0$  such that
$$
\sup_{x\in K, t\in [0,S]} E_xe^{\tilde \alpha\tau_K^t}<\infty, \quad \tau_K^t:=\inf \{s\geq 0:
X_{t+s}\in K\}.
$$
\end{itemize}

 In order to
simplify the exposition, we consider the case $\XX=\Re^d,$ only. One
can easily extend the proof to the general case by a standard
localization procedure.

We put
$\phi(x)=E_xe^{\tilde \alpha\tau_K}, \psi(x)= E_xe^{\alpha\tau_K}, x\in\XX$. Let us show that
$\phi$ is locally bounded; this would imply the condition (a) above.

Let $x_0\in \Re^d$ and $0<r_0<r_1$ be such that $K\subset
\{x:\|x-x_0\|<r_0\}$. Denote $D=\{x:\|x-x_0\|<r_1\}\setminus K$,
$\theta=\inf\{t:X_t\in \prt D\}$, and  $\mu_{x}(dy)\eqdef
P_x(X_\theta\in dy), x\in D$.

Consider auxiliary function
$$
h(x)=\int_{\prt D}E_{y}e^{\alpha \tau_K}\, \mu_x(dy)\in [0,\infty],\quad x\in D.
$$
This function can be represented as a monotonous point-wise limit of the functions
$$
h_N(x)=\int_{\prt D}g_N(y)\, \mu_x(dy), \quad N\geq 1
$$
with bounded and measurable functions $g_N$. Every function is $A$-harmonic in $D$, this  can be proved in a standard way using the strong Markov property of $X$, e.g. Chapter II \S5, \cite{DY}. Hence every $h_N$ satisfies the
Harnack inequality (see \cite{KS81}). Namely, there exists
$C\in\ax$ independent of $N$  such that
$$
h_N(x_1)\leq Ch_N(x_2)
$$
for every $y\in D$, and $x_1,x_2\in \{x: \|x-y\|<{1\over
2}\mathrm{dist}(y,\prt D)\}$. Then the same relation holds true with $h$ instead of $h_N$.  On the other hand, by the strong
Markov property of $X$, we have
$$
E_xe^{\alpha\tau_K}=E_x(e^{\alpha\theta}\psi(X_\theta))\geq
E_x\psi(X_\theta)=h(x), \quad x\in D.
$$
Hence, under condition $3')$, $h(x)<+\infty$ for $\pi$-a.a. $x\in
D$. In addition,  $\mathrm{supp}\, \pi=\XX$;  one can easily
verify this fact using positivity of the transition probability
density. Therefore, the function $h$ is bounded on every compact
$S\subset D$.

The function $h$ can be written in the form
$$
h(x)=E_xe^{\alpha\tau_K^\theta},\quad \tau_K^\theta=\inf\{s\geq 0:
X_{s+\theta}\in K\}.
$$
For $x\in D$, we have $\tau_K=\theta+ \tau_K^\theta$ $P_x$-a.s.,
and therefore
$$
E_xe^{\tilde \alpha\tau_K}\leq [E_x(e^{{\alpha\tilde\alpha\over
\alpha-\tilde\alpha}\theta})]^{\alpha-\tilde \alpha\over
\alpha}[h(x)]^{\tilde \alpha \over \alpha}.
$$
Using the Kac formula one can show that, for every $a>0$,  the
function $x\mapsto E_x e^{a\theta}$ is bounded on $D$ (this fact
is quite standard and hence we do not go into details here).
Therefore, the function $\phi$ is bounded on every compact
$S\subset D$.

Next, consider the closed ball $E=\{x:\|x-x_0\|\leq r_0\}$; note that its boundary  $S=\{x:\|x-x_0\|= r_0\}$ is a compact subset of $D$ and therefore the function $\phi$ is bounded on $S$. We put
$\sigma=\inf\{t:X_t\in S\}$, then by the strong Markov property of $X$ we have for $x\in E$
$$
\phi(x)\leq E_x(e^{\tilde \alpha\sigma}\phi(X_\sigma))\leq (E_x
e^{\tilde \alpha\sigma})\sup_{y\in S}\phi(y).
$$
The function $x\mapsto E_x e^{\tilde \alpha\sigma}$ is bounded on
$E$ (again, we do not give a detailed discussion here). Hence
$\phi$ is bounded on $E$. Since $r_0$ and $r_1$ can be taken arbitrarily
large, this means that  $\phi$ is locally bounded.

Now, let us  verify the condition (b) above.  We keep the notation $E=\{x:\|x-x_0\|\leq r_0\}, S=\prt E$, and put
$\sigma^0=0$,
$$
\sigma^{2n-1}=\inf\{t\geq \sigma^{2n-2}:X_t\in S\},\quad
\sigma^{2n}=\inf\{t\geq \sigma^{2n-1}:X_t\in K\},\quad n\geq 1.
$$
For any  $a>0$, one has   $$ q\eqdef \max\left[\sup_{x\in
K}E_xe^{-a\tau_S}<1, \sup_{x\in S}E_xe^{-a\tau_K}<1\right]<1
$$
because $\mathrm{dist}\,(K,S)>0$ and $X$ is a  Feller process
with continuous trajectories. Therefore,
 \be\label{83}
E\Big[e^{-a(\sigma^{k+1}-\sigma^k)}\Big|\Ff_{\sigma^k}\Big]\leq
q\quad \hbox{a.s.,} \quad k\geq 0. \ee We have
$$
E_xe^{\tilde \alpha\tau_K^t}=\sum_{k=0}^\infty E_xe^{\tilde
\alpha\tau_K^t}\1_{\sigma^k\leq t<\sigma^{k+1}}, \quad x\in K.
$$
For $k$ even,  $X_t\in E$ a.s. on the set
$C_{k,t}\eqdef\{\sigma^k\leq t<\sigma^{k+1}\}$. In addition,
$C_{k,t}\in \Ff_t$. Hence
$$\ba
E_xe^{\tilde \alpha\tau_K^t}\1_{\sigma^k\leq t<\sigma^{k+1}}&=
E_x\left(\1_{\sigma^k\leq t<\sigma^{k+1}}E\Big[e^{\tilde
\alpha\tau_K^t}\Big|\Ff_t\Big]\right)= E_x\1_{\sigma^k\leq
t<\sigma^{k+1}}\phi(X_t)\\
&\leq \sup_{y\in E}\phi(y)\,P_x(\sigma^k\leq t<\sigma^{k+1}),
\quad k=2n. \ea
$$
For $k$ odd,  $\tau_K^t=\sigma^{k+1}-t\leq \sigma^{k+1}-\sigma^k$
a.s. on the set $C_{k,t}$. Hence
$$\ba E_xe^{\tilde \alpha\tau_K^t}\1_{\sigma^k\leq t<\sigma^{k+1}}&\leq
E_x\1_{\sigma^k\geq t}e^{\tilde
\alpha(\sigma^{k+1}-\sigma^k)}=E_x\left(\1_{\sigma^k\leq
t}E\Big[e^{\tilde
\alpha(\sigma^{k+1}-\sigma^k)}\Big|\Ff_{\sigma^k}\Big]\right)\\
&=E_x\1_{\sigma^k\leq t}\phi(X_{\sigma^k})\leq \sup_{y\in
E}\phi(y)\, P_x(\sigma^k\leq t).\ea
$$
Therefore,
$$
E_xe^{\tilde \alpha\tau_K^t}\leq \sup_{y\in
E}\phi(y)\sum_{k=0}^\infty P_x(\sigma^k\leq t), \quad x\in K.
$$
It follows from (\ref{83}) that $E_x e^{-a\sigma^k}\leq q^k, x\in
K$. Then
$$
P_x(\sigma^k\leq t)=P_x(-\sigma^k\geq -t)\leq e^{at}q^k, \quad
k\geq 0, x\in K,
$$
and consequently
$$
\sup_{x\in K, t\in [0,S]} E_xe^{\tilde \alpha\tau_K^t}\leq
e^{aS}(1-q)^{-1}\sup_{y\in E}\phi(y)<+\infty.
$$
We have verified  conditions (a) and (b).  Hence the required statement  follows from  Proposition
2.4 in \cite{Kul11}.\qed

\end{document}